\documentclass[11pt]{article}

\usepackage{amsmath, amscd}
\usepackage{amssymb}
\usepackage{diagrams}
\newcommand{\tto}{\rightrightarrows}

\newtheorem{prop}{Proposition}

\newtheorem{lemma}[prop]{Lemma}
\newtheorem{thm}[prop]{Theorem}
\newtheorem{definition}[prop]{Definition}

\newcommand{\pg}{parallelogram}
\newcommand{\pp}{par\-al\-lel\-e\-pi\-pe\-dum}
\newcommand{\ppa}{par\-al\-lel\-e\-pi\-peda}

\DeclareMathOperator{\sign}{sign}

\DeclareMathOperator{\proj}{proj}
\DeclareMathOperator{\id}{id}
\DeclareMathOperator{\Vol}{Vol}

\begin{document}

\title{Infinitesimal cubical structure, and higher connections}
\author{Anders Kock}
\date{}

\maketitle
\section*{Introduction} The purpose of the present note is to 
experiment with 
a possible framework for the theory of ``higher connections'', as 
have recently become expedient in string theory; however, our work 
was not so much to find a framework which fits any existing such 
theory, but rather to find notions which come most naturally of their 
own accord. The approach is based on the combinatorics 
of the ``first neighbourhood of the diagonal'' of a manifold, using 
the technique and language of Synthetic Differential Geometry (SDG), 
as in \cite{SDG}, and notably in \cite{CCBI}. The present note may be 
seen as a sequel to the latter, and also to \cite{CPCG}.
 The basic viewpoint in \cite{CCBI} is 
that connections (1-connections) take value in {\em groupoids} (a 
viewpoint which goes back to Ehresmann), and that they in effect may 
be seen as morphisms in the category of reflexive symmetric graphs, 
noting that any groupoid has an underlying such.

To go beyond this into higher dimensions, one would like to consider 
some kind of higher groupoids to receive the values of the 
connections (see e.g.\ \cite{schreiber}), and 
we also need a higher dimensional version of the notion of reflexive 
symmetric graph. This led us to pass into the cubical, rather than into the 
simplicial world. The passage into this world depends on having a 
cubical complex associated to any manifold, in analogy with the 
simplicial complex of ``infinitesimal simplices'' that one 
derives out of the first neighbourhood of the diagonal. This latter 
simplicial complex is known to be the carrier of a theory of 
``combinatorial differential forms'', as in  \cite{SDG}, \cite{CCBI}, 
\cite{DFIC}, and in  \cite{BM}.

The observation that opens up for a similar cubical complex is that 
infinitesimal simplices in a manifold canonically give rise to 
infinitesimal \ppa . This hinges on the possibility of forming affine 
combinations of the points in an infinitesimal simplex, a possibility 
first noted in \cite{GCLCP}, see also \cite{DFIC}.

On the algebraic side, the kind of ``higher'' groupoid which fits the 
bill are essentially the ``$\omega$-groupoids'' of Brown and Higgins 
\cite{BH}, or their truncation to some finite dimension, in which 
case we call them $n$-cubical groupoids. For $n=2$ they are the 
edge-symmetric double groupoids of Brown and Spencer, \cite{BS}. In 
any case, all structure involved (at present) is strict, and no coherence issues 
are involved. 

Besides connection as an infinitesimal notion, we study a 
corresponding global notion, which is that of holonomy. We relate 
holonomy of $n$-connections with integrals of differential $n$-forms.

This research was partly triggered by some questions which 
Urs Schreiber posed me in 2005; for $n=1$, an attempt of an answer 
was provided in my \cite{CPCG}. I want to thank him for the impetus. I 
also want to thank Ronnie Brown for having for many years persuaded 
me to think strictly and cubically. Finally, I want to thank Marco 
Grandis for useful conversations on cubical and other issues.

\section{Preliminaries on cubical sets}\label{preliminaries}
The fundamental role of cubical sets (=cubical complexes) have been 
emphasized  through the work of Brown and his collaborators, 
and by Grandis and his collaborators, and by many others,
 but is generally not so well 
documented in the literature. So we start with a brief ``lexicon'' to 
establish terminology. This in particular applies to the  
``symmetry'' structures which a cubical set may have, and to the 
extra ``degeneracy''-structure, called ``connections'' by Brown, 
Spencer and Higgins. Since the present paper deals with connections 
in a quite different sense of the word, we shall here use the term 
``BSH-connections'' for the kind of connections which they consider.

We shall mainly follow the scheme of terminology and notations as in 
\cite{GM}. Recall that a cubical set $C_{\bullet}$ is a family of 
sets $C_{n}$ (the ``set of {\em $n$-cubes} or {\em $n$-cells}'';
 $n$ ranging over the natural numbers), equipped with 
face- and degeneracy operators:
$$\mbox{face operators } \partial^{\alpha}_{i}:C_{n}\to C_{n-1} \quad 
(i=1, \ldots ,n, \alpha = 0,1)$$
$$\mbox{degeneracy operators } \epsilon _{i}:C_{n-1}\to C_{n} \quad 
(i=1, \ldots ,n)$$
satisfying the ``cubical relation''-equations, as in \cite{BH} (1.1) 
(or \cite{GM} (5) for the dual relations), in analogy with the 
standard relations for the face- and degeneracy operators for 
simplicial sets.

Just as a simplicial set may carry a richer structure, namely  
compatible 
actions of the symmetric groups ${\mathfrak S}_{n+1}$ on the set of 
$n$-simplices (in which case it is called a {\em symmetric} simplicial 
set by \cite{GTAC8}), a cubical set may carry certain symmetry 
structures; the symmetries now come in two classes, called  
respectively, {\em 
interchanges} or {\em transpositions}, and {\em reversions} in 
\cite{GM}.

 The {\em transpositions} 
which a cubical set may have give rise to an action of ${\mathfrak 
S}_{n}$ on the set $C_{n}$. This action may be given in terms of 
operators $\sigma _{i}: C_{n}\to C_{n}$ ($1\leq i\leq n-1$) 
satisfying certain relations and compatibilities with the face- and 
degeneracy operators, (see \cite{GM}, equations (27)-(29)). 
Geometrically, $\sigma _{i}$, applied to an $n$-cube, interchanges the 
$i$th and $(i+1)$st edge emanating from the initial vertex.

The {\em reversions} are given in terms of operators $\rho 
_{i}:C_{n}\to C_{n}$ ($1\leq i\leq n$) satisfying 
certain relations, and compatibilities with the face- and 
degeneracy operators; and also compatibilities with the interchange 
operators, in case this structure is present, (see \cite{GM}, 
equations (56)). Geometrically, $\rho _{i}$ applied to an $n$-cube, 
performs a reflection of the cube in the hyperplane orthogonal to the 
$i$th axis.

(In case a cubical set has both transpositions and reversions in a 
compatible way, the ``hyperoctahedral group $({\mathbb 
Z}/2)^{n}\rtimes {\mathfrak S}_{n}$'' acts on the set of $n$-cubes, 
see \cite{GM} \S 9.)

In the theory of connections that we develop, the reversions are 
more important than the transpositions. For, in the cubical groupoids 
that we consider below, $\rho _{i}$ is canonically present in terms of the 
{\em inversion} of arrows.

\medskip

Finally, a cubical set may have BSH-connections; these are an extra 
family of degeneracy operators $\gamma _{i}:C_{n-1}\to C_{n}$ ($1\leq 
i \leq n-1$), with compatibility relations with the other structural 
elements (face-, degeneracy, transposition-, and/or reversion 
operators). They are denoted $\Gamma _{i}$ in \cite {BH}, and we 
give some comments and geometry about these in Section \ref{BSHT} below.
The cubical set of ``infinitesimal \ppa '' arising from a manifold, 
which is the main carrier of the theory we present here, has
\mbox{face-,} degeneracy-, transposition-, and reversion- 
operators, but does not have BSH-connections. (The same applies to the 
cubical set of  \ppa\ in an arbitrary affine space. There is also a 
notion of \ppa\ in a non-commutative affine space (pregroup); here, neither 
BSH-connections nor transpositions are present. This digression, we 
include in Section \ref{digression} below.)

\subsection{Shells}
Let $Cub$ denotes the category of cubical sets, and $Cub _{n}$ the 
category of cubical sets ``up to dimension $n$'' (such things, we 
call $n$-cubical sets). There is an evident 
truncation functor $tr: Cub _{n+1}\to Cub _{n}$, forgetting the 
$n+1$-cells. It is a functor between two presheaf categories, and is 
induced by a functor between the respective index-categories, and it 
therefore has adjoint on both sides. The right adjoint (=the 
coskeleton), we denote 
just by $(-)'$ (following \cite{BH}), thus if $H$ is an $n$-cubical 
set, $H'$ is an $n+1$-cubical set; for $k\leq n$, $H' _{k}= H_{k}$, 
whereas $H' _{n+1}$ consists of ``$n$-shells'' ${\bf x}$ in $H_{n}$, meaning 
$2(n+1)$-tuples of elements in $H_{n}$, whose boundaries match up as if these
 cells were the $2(n+1)$ faces of an $n+1$-cell in some cubical set. 
These $n$-cells then serve as the faces of ${\bf x}$ in $H'$.

If $K$ is an $n+1$-cubical set, there is a front adjunction for the 
adjointness considered,
$\theta : K\to (tr(K))'$; it is the identity in dimensions $\leq n$, 
and in dimension $n+1$, it associates to an $n+1$-cell $x$ in $K$ its 
boundary shell ${\bf \partial}(x)$. (These ideas are from 
\cite{BH}.)

All this lifts to cubical sets with reversions and/or transpositions 
and/or BSH connections. 
We shall not introduce special notation for 
each of these enriched notions, but indicate in each case which one 
of these notions we consider. 

\subsection{The complex of singular cubes}\label{complex}

If $M$ is a manifold, we define the set $S_{k}(M)$ of singular 
$k$-cubes in $M$ to be the set of (smooth) maps $f:R^{k}\to M$. (The reader may 
prefer to think of the relevant information of the singular cube $f$ 
to reside in the restriction of $f$ to the unit cube $I^{k}\subseteq 
R^{k}$, but this 
unit cube does not play a role in the formalism we present.)
The $S_{k}(M)$ 
jointly carry 
structure of a cubical complex $S_{[\bullet]}(M)$: for $\alpha = 0$ or $=1$, and 
$i= 1, \ldots , k$, the face map $\partial ^{\alpha}_{i} : S_{k}(M)\to 
S_{k-1}(M)$ consists in precomposing with the affine map $R^{k-1}\to R^{k}$, 
$(x_{1}, \ldots ,x_{k-1})\mapsto (x_{1}, \ldots , \alpha , \ldots 
,x_{k-1})$ ($\alpha = 0$ or $=1$, placed in the $i$th position),
 and the degeneracy map $\epsilon 
_{i}:S_{k-1}(M)\to S_{k}(M)$ is induced by  projection $R^{k}\to 
R^{k-1}$ (omitting $i$th coordinate). This is also an affine map. 
There are further canonical affine, hence smooth,  maps between the 
$R^{n}$s, giving a richer strucure to $S_{[\bullet]}(M)$, namely 
transposition of the coordinates, and reversions.
 We shall also consider these. But one kind of 
structure will not be present in the purely affine world  to be 
considered next, namely the (BSH-) connections
$\gamma _{i}:S_{k-1}(M)\to S_{k}(M)$.

\medskip

If $M$ is an {\em affine} space, there is a subcomplex $A_{[\bullet 
]}(M)$ of 
$S_{[\bullet]}(M)$, with $A_{[n]}(M)$ consisting of the {\em affine} maps
$R^{n}\to M$ (``{\em affine} singular $n$-cubes''). Reversions and transpositions of $S_{[\bullet]}(M)$ 
restrict to this subcomplex. The information of an affine singular 
$n$-cube is of course contained in the images of the $2^{n}$ vertices 
of the unit cube, and they form the vertices of an $n$-dimensional 
\pp\ in $M$.

\medskip

By an $n$-{\em simplex} in a set $M$, we understand an $n+1$-tuple of 
points in $M$, called the {\em vertices} of the simplex.
Since $R^{n}$ is a free affine space on the $n+1$-tuple of points $0, e_{1 }, 
\ldots ,e_{n}$, it follows that for an {\em affine} space $M$, there is a 1-1 
correspondence between $n$-simplices in $M$, and affine singular 
$n$-cubes. We need a notation: to an $n$-simplex $(x_{0}, \ldots ,x_{n})$ 
in $M$, we denote the corresponding affine map $R^{n}\to M$ by 
$[[x_{0}, \ldots ,x_{n}]]$; it is given explicitly by the formula
\begin{equation}[[x_{0}, \ldots ,x_{n}]](s_{1}, \ldots ,s_{n})= (1-\sum_{1}^{n} s_{i})\cdot x_{0} + 
\sum_{1}^{n}s_{i}\cdot x_{i}.\label{aff}\end{equation}
Note that the right hand side here is an affine combination, i.e, the sum of the 
coefficients is 1.
Thus, the graded set of simplices in an affine space carries 
the structure of a (symmetric) simplicial set, as well as of a cubical 
set (with reversions and transpositions).

\begin{sloppypar}
Another way to encode the information of an affine singular $n$-cube
$[[x_{0}, \ldots ,x_{n}]]$ 
in an affine space $M$ is to consider the $2^{n}$-tuple of points 
whih are the images of the vertices of the unit cube in $R^{n}$
under the map $[[x_{0}, \ldots ,x_{n}]]: R^{n}\to M$. They form the 
vertices of an $n$-dimensional \pp\ in $M$. (This viewpoint will be 
elaborated algebraically in the Digression below.) The $2^{n}$-tuple of these points we denote $P(x_{0};x_{1}, 
\ldots ,x_{n})$; it contains exactly the same information as the 
simplex $(x_{0}, \ldots ,x_{n})$, since the vertices of the latter appear in a 
canonical way as (some of) the elements of $P(x_{0};x_{1}, 
\ldots ,x_{n})$. The $2^{n}$-tuple $P(x_{0};x_{1}, 
\ldots ,x_{n})$, we call an {\em $n$-dimensional \pp\ } in $M$.
\end{sloppypar}

\medskip

We shall elaborate a little of properties of the bijection between 
simplices and \ppa\ in an affine space $M$, because 
these properties will also apply for the {\em infinitesimal} simplices and 
\ppa\ which we shall consider. So let $M$ denote an affine space; 
let $s_{(\bullet )}(M)$ denote the symmetric simplicial set of simplices 
in it, and 
$A_{[\bullet ]}(M)$ the cubical set, with  transpositons and 
reversions, consisting in the \ppa\ in $M$. Then we have for each $n$ 
the comparison map $p_{n}: s_{(n)}(M)\to A_{[n]}(M)$ defined by
$$(x_{0},x_{1},\ldots ,x_{n})\mapsto P(x_{0};x_{1}, \ldots ,x_{n}).$$
The group ${\mathfrak S}_{n+1}$ acts on $s_{(n)}(M)$, and ${\mathfrak 
S}_{n}$ acts as the subgroup of permutations which fix the vertex 
$x_{0}$.
Then it is clear that the comparison map $p_{n}:s_{(n)}(M)\to 
A_{[n]}(M)$ preserves
 the action of ${\mathfrak S}_{n}$.
The question of degenerate simplices vs.\ degenerate cubes is more 
subtle. A simplex with $x_{i}=x_{i+1}$ is degenerate in the sense of 
being a value of a degeneracy operator; whereas a \pp\ is a value of 
a degeneracy operator if $x_{0}=x_{i}$ for some $i>0$. The degenerate 
2-simplex $(x_{0},x_{1},x_{1})$ goes by $p$ to a \pg\ which is not a 
value of a cubical degeneracy operator, not even modulo some action of the 
symmetry group. In fact, for $x_{0}\neq x_{1}$ it has {\em three}
 distinct vertices, namely $x_{0}$, 
$x_{1}$ and $2x_{1}-x_{0}$.

\subsection{Subdivision}\label{subdivision}
We return to the general case where $M$ is a manifold, not 
necessarily an affine space.
We shall  consider a structure which $S_{[\bullet]}(M)$ has (we 
don't here attempt to axiomatize this structure); it is the notion of 
{\em subdivision}. Namely, given a scalar $s\in R$, 
and an index $i=1, \ldots ,n$, we have a ``subdivision'' operation, 
which to a singular $n$-cube $f:R^{n }\to M$ associates a pair of 
singular $n$-cubes $f'$ and $f''$; they come about by precomposing 
$f$ by certain affine maps $h_{s,i}:R^{n}\to R^{n}$ and $k_{s,i}:R^{n}\to 
R^{n}$, respectively. Here $h_{s,i}$ is the affine map corresponding to the 
$n$-simplex $(0, e_{1}, \ldots ,s\cdot e_{i}, \ldots ,e_{n})$, and 
$k_{s,i}$ is the affine map corresponding to the $n$-simplex $(s\cdot e_{i},e_{1}+s\cdot 
e_{i}, \ldots ,e_{i}, \ldots ,e_{n}+s\cdot e_{i})$. Alternatively, 
$h_{s,i}:R^{n}\to R^{n}$ is the the map $\id \times \ldots \times 
h_{s}\times \ldots \times \id $, where $h_{s}:R\to R$ is the map 
$t\mapsto t\cdot s$; similarly for $k_{s,i}$ with $k_{s}:R\to R$ the 
map $t\mapsto s + t\cdot (1-s)$.

Let us make a picture, and thereby also introduce a notation which we 
shall need later on. We consider the case of an affine space $M$ (we 
shall ultimately be interested in the plane $R^{2}$), and consider an 
{\em affine} singular 2-cube $[[x,y,z]]$, geometrically a \pg\ with 
vertices $x,y,z,u$.
We consider for each $s\in R$ the point 
$y(s):=(1-s)x +sy$, and also the point $u(s):=(1-s)z + su$. 

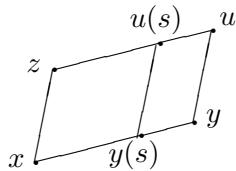
\begin{equation}
\begin{picture}(100,75)(0,-10)
\put(20,6){\line(4,1){60}}
\put(20,6){\line(1,5){7}}
\put(80,21){\line(1,5){7}}
\put(27,41){\line(4,1){60}}
\put(10,4){$x$}
\put(85,21){$y$}
\put(17,40){$z$}
\put(19,5){{\bf .}}
\put(79,20){{\bf .}}
\put(26,40){{\bf .}}
\put(86,55){{\bf .}}
\put(90,58){$u$}
\put(59,15){{\bf .}}
\put(59,15){\line(1,5){7}}
\put(48,6){{$y(s)$}}
\put(66,50){{\bf .}}
\put(56,57){$u(s)$}
\end{picture}
\label{subdiv}\end{equation}

\noindent For each 
$s,t$, $P(x;y(s+t),z)$ subdivides into $P(x;y(s),z)$, $P(y(s); 
y(s+t),u(s))$. In particular, the displayed \pg\ 
$P=P(x;y,z)$ subdivides into
$P'=P(x;y(s),z)$ and $P'' = P(y(s);y,u(s))$.

If $f'$ and $f''$ come about from $f$ by such subdivision process (with $s$ 
and $i$ fixed), we say that $f$ $(i,s)$-{\em subdivides into} $f'$, 
$f''$. In this case
we have the following equations 
\begin{equation}\partial ^{1}_{i}(f')=\partial 
^{0}_{i}(f'')\label{subd1}\end{equation}
\begin{equation}\partial ^{0}_{i}(f')=\partial 
^{0}_{i}(f)\label{subd2}\end{equation}
\begin{equation}\partial ^{1}_{i}(f'')=\partial 
^{1}_{i}(f).\label{subd3}\end{equation}
Note that these equations are exactly the book-keeping conditions 
that one has if $f$ is the $i$-composite of $f'$ and $f''$ in a 
cubical groupoid (in the sense recalled in the next Section).
There are also 
 equations relating the $j$-faces of $f$ with the faces of 
$f'$ and $f''$: the formulae are, for $\alpha =0$, and for 
$\alpha =1$,
\begin{equation}(\partial ^{\alpha}_{j}(f))' =\partial ^{\alpha} _{j}(f')
\label{subd4}\end{equation}
\begin{equation}(\partial ^{\alpha}_{j}(f))'' =\partial ^{\alpha} _{j}(f'')
\label{subd5}\end{equation}
and similarly for $f''$; but for $j<i$, the prime  on the left hand 
side of (\ref{subd4}) refers to $h_{s,i-1}$, and the double prime 
on the left hand side of (\ref{subd5})similarly refers to $k_{s,i-1}$. 
These equations are all seen simply by considering the respective 
representing (affine) maps between the relevant $R^{m}$s.

\medskip

These equations can also be formulated
\begin{prop}Assume $f$ $(i,s)$-subdivides into $f'$, $f''$. Then for
$j>i$, $\partial ^{\alpha}_{j}(f)$ $(i,s)$-subdivides into $\partial 
^{\alpha}_{j}(f')$, $\partial ^{\alpha}_{j}(f'')$, and for $j<i$,
$\partial ^{\alpha}_{j}(f)$ $(i-1,s)$-subdivides into $\partial 
^{\alpha}_{j}(f')$, $\partial ^{\alpha}_{j}(f'')$
\label{subdivisionp}\end{prop}

\subsection{Cubical groupoids}
We take the notion of {\em cubical groupoid} from the fundamental paper 
\cite{BH} by Brown and 
Higgins. It is what they call $\omega$-groupoids, except that we do 
not assume that their ``connections'' $\Gamma _{j}$ are part of the 
structure. The term ``cubical groupoid'', we have from Grandis.
 Thus, a cubical groupoid $G$ is a cubical complex, where 
each $G_{n}$ is equipped with $n$ (partially defined) compositions 
$+_{i}$; the composition $+_{i}$ makes $G_{n}$ into the set of arrows 
of a groupoid, with $G_{n-1}$ as set of objects and with
 $\partial ^{0}_{i}$ and $\partial ^{1}_{i} :G_{n}\to G_{n-1}$ as 
domain- and codomain-formation, respectively. The degeneracy map 
$\epsilon _{i}:G_{n-1}\to G_{n}$ provides the identity maps for this 
groupoid structure. Formation of inverse arrows w.r.to $+_{i}$ is 
denoted $-_{i}$ and makes $G$ into a cubical set with reversions. 
There are several compatibility equations between the structural 
elements; they may be read off in \cite{BH} 1.3, 1.4, and 1.5.

The notion of cubical groupoid truncates in an evident way:
thus, by an {\em $n$-cubical groupoid}, we understand an $n$-cubical set 
$G$ with compositions $+_{i}$ ($i=1,\ldots ,k$) on each $G_{k}$ 
($k\leq n$), satisfying the same family of compatibility equations.
1-cubical groupoids are just groupoids;
2-cubical groupoids are double groupoids where the set of horizontal 
and vertical 1-cells are equal (the ``edge-symmetric'' double 
groupoids of Brown and Spencer \cite{BS}). The truncation functors $tr: 
Cub_{n+1}\to Cub _{n}$ and their right adjoints $(-)' :Cub _{n}\to Cub _{n+1}$ 
 lift to functors between the categories of $n$- and 
$n+1$-cubical groupoids. For $tr$, this is trivial; for $(-)'$: if 
$G$ is an $n$-cubical groupoid, $G'_{n+1}$ consists of {\em shells} of 
cells in $G_{n}$, and they may be concatenated, or reversed, in each of the $n+1$ 
directions, using the compositions and reversions of $G_{n}$. The 
exact formulae may be found in \cite{BH} \S 5. In particular, for 
$n=1$: if $G$ is an ordinary groupoid, $G'$ is the familiar (edge 
symmetric) double groupoid 
whose 2-cells are the (not necessarily commutative) squares in $G$.

\medskip

The affine singular cubes in an affine space $A$ form a cubical groupoid; 
we shall describe its  $i$th composition. Given $n$-cubes
(i.e.\ $n$-dimensional \ppa ) in $A$, say $P$ and $P'$ 
with $\partial _{i}^{1}(P) =\partial _{i}^{0}(P')$, we 
describe the 
$n$-cube $P\circ _{i}P'$ as follows.
Let $P=P(x_{0};x_{1}, \ldots ,x_{n})$ and similarly
$P'=P(x_{0}';x_{1}', \ldots ,x_{n}')$. The compatibility
$\partial _{i}^{1}(P) =\partial _{i}^{0}(P')$ forces, for 
$j\neq i$, the 
$x_{j}'$ to equal certain combinations of the $x_{k}$s. Then
$$P\circ _{i}P'=(x_{0}; x_{1}, \ldots ,x_{i}', \ldots 
,x_{n}).$$
Note that if $P$ is $(i,s)$-subdivided into $P'$, $P''$, then $P=P' 
+_{i}P''$.
-- Degenerate cubes $s_{i}(c)$ act as identities for $\circ _{i}$. We 
finally describe inversion $-_{i}$ in the $i$th direction:
$$-_{i}P(x_{0};x_{1}, \ldots ,x_{n})= P(x_{i};x_{1}-x_{0}+x_{i}, \ldots 
,x_{0}, \ldots , x_{i}-x_{0}+x_{n})$$
with the $x_{0}$ placed in the $i$th slot.

It is easy to verify that this decribes a cubical groupoid (in the 
sense of \cite{BH}). (It is in fact a cubical equivalence relation: 
for any two $n-1$-cubes and any $i=1, \ldots n$, there is at most 
one $n$-cube $c$ having the given cubes 
as $\partial _{i}^{0}$- and $\partial _{i}^{1}$-faces, respectively.)

This cubical groupoid also has transpositions: $\sigma 
_{1}(P(x_{0};x_{1}, x_{2},\ldots ,x_{n}) = P(x_{0};x_{2},x_{1}, 
\ldots ,x_{n})$ etc. They will not exist in the non-commutative case 
which we now consider, as a digression:

\subsection{Digression on the non-commutative case}\label{digression}
The cubical groupoid structure described above for an affine space $A$ 
also exists for ``non-commutative affine spaces'', i.e.\ for 
``pregroups'' in the sense of \cite{ATMF}; to make the exposition 
more accessible to the reader not aquainted with pregroups,
 we consider the special case of a group instead of a pregroup.
So let us consider a group $G$, not necessarily commutative, but we use 
additive notation. To an $n+1$-tuple $x_{0}, x_{1}, \ldots ,x_{n}$ in 
$G$, we associate a $2^{n}$-tuple of elements of $G$ in the following 
way: let $h_{1}<h_{2}<\ldots <h_{k}$ be a $k$-element subset $H$  of 
the set $\{ 1, \ldots ,n\} $. If $k\geq 2$, we associate to it the element 
$x_{H}$
of $G$ 
given by the expression
$$x_{h_{1}}-x_{0}+x_{h_{2}}-x_{0}+\ldots -x_{0}+x_{h_{k}};$$
if $H$ is a singleton subset $H=\{ h \}$, we put $x_{H}=x_{h}$, and 
if $H=\emptyset $, we put $x_{H}=x_{0}$. Note that in the 
commutative case, all these expressions are affine combinations in 
$G$, in fact, the $2^{n}$-tuple descibed is exactly the $2^{n}$-tuple 
of vertices in the $n$-dimensional \pp\ $P(x_{0};x_{1},\ldots ,x_{n})$. 
Then 
$$\partial ^{0}_{i}(P(x_{0};x_{1}, \ldots ,x_{n}))=P(x_{0}; x_{1}, 
\ldots ,\widehat{x_{i}}, \ldots ,x_{n})$$
whereas
$$\partial ^{1}_{i}(P(x_{0};x_{1}, \ldots ,x_{n})) =P(x_{i}; x_{\{ 
1,i \} }, 
\ldots ,\widehat{i}, \ldots ,x_{\{ i,n \} }).$$
Degeneracies consist in inserting $x_{0 }$ ($i\geq 1$); the formulae 
for composition and reversion are the same as described above for 
affine space.
-- The transpositions described for the case of affine space $A$ will 
not work here. For instance, for 
$n=2$, there is one transposition operator $\sigma$, and $\partial 
^{1}_{1}\circ \sigma = \partial ^{1}_{2}$ is one of the required 
compatibilities. Consider $P(x;y,z)$. Then
$$\partial 
^{1}_{1}( \sigma (P(x;y,z)))=\partial ^{1}_{1}P(x;z,y) = P(z;z-x+y)$$
whereas
$$\partial ^{1}_{2}(P(x;y,z))=P(z;y-x+z).$$
These may be different, in the non-commutative case.

\subsection{Brown-Spencer-Higgins theory}\label{BSHT}
If $C=C_{\bullet}$ is a cubical set, there is a graded subset $Cr 
(C)$ of $C$, where $Cr _{n}(C)\subset C_{n}$ for $n\geq 2$ consists of those $c$ 
such that all faces of $c$ except possibly $\partial ^{0}_{1}(c)$ are 
totally degenerate, i.e.\ come about by applying $n$ degeneracy 
operators to a $0$-cell $\in C_{0}$. Then $\partial ^{0}_{1}$ 
restricts to a map $\delta : Cr_{n}(C)\to Cr _{n-1}(C)$. If $C$ is 
equipped with structure of a cubical groupoid, the graded set $Cr 
_{\bullet}(C)$ carries the structure of a {\em crossed complex} over 
the groupoid $C_{1}\tto C_{0}$ in the sense of \cite{BH}, in 
particular, each $C_{n}$ for $n\geq 2$ is a $C_{0}$-indexed family of 
groups (abelian, if $n\geq 3$), and $\delta$ is a group homomorphism, 
with $\delta \circ \delta =0$.

We need to recall the theory of connections in the sense of 
Brown, Spencer, and Higgins. We  
call them {\em BSH-connections}. They are extra degeneracies $\gamma 
_{i}: C_{n-1}\to C_{n}$ ($i=1, \ldots ,n-1$) with which a 
cubical set $C_{\bullet}$ may be equipped. For instance, the 
2-cubical groupoid $G'$ of squares in an ordinary 1-groupoid $G$ has 
the extra degeneracy operator $\gamma$ which to an arrow $g:x\to y$ 
in $G$ 
associates the square in $G$
$$\begin{diagram}x&\rTo ^{g}&y\\
\dTo ^{g}&&\dTo_{\id}\\
y&\rTo _{\id}&y.
\end{diagram}$$

If $G=G_{\bullet}$ is a cubical groupoid with BSH-connections (compatible with the 
groupoid structures, cf.\ \cite{BH}), there is a retraction 
(``folding operation'', \cite{BH}),   
 ${\bf \phi}: G_{n}\to Cr_{n}(G)$ (for all $n$).
The folding operations give a way of 
expressing what it means for an $n$-cell in a cubical groupoid to 
commute:

 If $G$ is an $n$-cubical groupoid with 
BSH-connections, we are interested in the question of commutativity of 
$n+1$-cells  in the $n+1$-cubical groupoid $G'$; 
recall that $G'$ has BSH-connections if $G$ does, so that it has its own 
folding operations. To say that ${\bf x}\in G'_{n+1}$ is commutative 
is then taken to mean 
that ${\bf \phi}({\bf x})$ is the zero element of (one of the  
groups which constitute) $Cr_{n+1}(G')$. Recall that an element of
 $G'_{n+1}$ is a shell of $n$-cells from $G$. E.g.\ for $n=1$, it is 
a square of arrows from the groupoid $G$, and the commutativity in 
the Brown-Higgins sense says that the cyclic composite
of the four arrows or their inverses is an identity arrow, and this 
property 
is equivalent with what everybody understands as a 
commutative square in a groupoid. 

\medskip

From Proposition 5.4 in \cite{BH}, we get
\begin{prop}Let $G$ be an $n$-cubical groupoid with BSH-connections;
then $Cr_{n+2}(G'') =\{ 0 \}$ (more precisely, it is the $G_{0}$-indexed 
family of 0-groups).
\label{BH5.4}\end{prop}

More information about the relationship between cubical groupoids and 
crossed complexes is considered in Section \ref{connections} below.

\section{Infinitesimal \ppa}\label{infinitesimal}We now place ourselves 
 in the context of Synthetic Differential Geometry. If $M$ is a 
manifold,  we have the notion of points $x,y$ in $M$ being (1st order) 
{\em neighbours}, written $x\sim y$. The neighbour relation $\sim$ is a reflexive symmetric 
relation. It is not transitive. A $k$-simplex $(x_{0}, x_{1}, \ldots 
,x_{k})$ in $M$ is called an {\em infinitesimal} simplex if 
$x_{i}\sim x_{j}$ for all $i,j = 0, \ldots ,k$.

It was 
proved in \cite{GCLCP} (see also \cite{DFIC}) that if $(x_{0}, x_{1}, \ldots 
,x_{k})$ is such an infinitesimal $k$-simplex, then affine 
combinations $\sum _{0}^{k} t_{i}\cdot x_{i}$ may be formed, using a 
coordinate chart, but independent of the chart; furthermore, any two 
of these affine combinations are neighbours. Also, any map $M\to N$ 
preserves the formation of such affine combinations. This implies that 
the formula (\ref{aff}) describing the map $[[x_{0}, \ldots ,x_{k}]]$, in case $M$ 
is an affine space, makes sense even if $M$ is just a manifold, provided 
$(x_{0}, \ldots ,x_{k})$ is an {\em ininitesimal} simplex. Thus, if
$(x_{0}, \ldots ,x_{k})$ is an infinitesimal simplex in a manifold 
$M$, it defines a singular $k$-cube $[[x_{0}, \ldots ,x_{k}]]$ in 
$M$. The singular $k$-cubes in $M$ which arise this way from 
infinitesimal simplices, we call {\em infinitesimal $k$-dimensional \ppa }.
The set of these is denoted $M_{[k]}$. The 
infinitesimal \ppa\ in $M$ form a subcomplex of the cubical complex 
$S_{[\bullet]}(M)$, stable under the reversion and transposition 
operations. Also, a subdivision of an infinitesimal \pp\ is an 
infinitesimal \pp .
 -- We denote this cubical complex $M_{[\bullet ]}$, and 
its truncation to $n$ dimensions, we denote $M_{[[n]]}$; the 
inclusion we denote $i_{n}$,
\begin{equation}\begin{diagram}M_{[[n]]}&\rTo ^{i_{n}}& S_{[[n]]}(M).
\end{diagram}\label{in} \end{equation}

The best geometric intuition of the cubical complex $M_{[\bullet ]}$ is that 
its $k$-cells are, or describe,  {\em infinitesimal $k$-dimensional \ppa\ }
in the following more geometric sense: if 
$(x_{0}, \ldots ,x_{k})$ is an infinitesimal $k$-simplex in $M$,  it 
defines a $2^{k}$-tuple of points, namely the images under $[[x_{0}, 
\ldots ,x_{k}]]:R^{k}\to M$ of the $2^{k}$ corner points of the unit 
cube $[0,1]^{k}$. To have this intuition playing an active role in 
our reasoning, we write $P(x_{0};x_{1}, \ldots ,x_{k})$ for this 
$2k$-tuple, even though it contains exactly the same information as 
the infinitesimal $k$-simplex $(x_{0}, \ldots ,x_{k})$, or as the singular $k$-cube 
$[[x_{0}, \ldots ,x_{k}]]$. But note that the ``unit interval'' 
$I=[0,1]$, as a point set $\{ s\in R \mid 0\leq s \leq 1 \}$, and 
similarly, the ``unit cube'' $I^{k}$, do not play any formal role in our 
treatment. (They can be incorporated in a more refined theory, where 
one takes a preordering $\leq$ of the line $R$ into account. If we 
let $\leq$ be the ``chaotic'' preordering ($x\leq y$ for all $x,y \in 
R$), then $R = [0,1]$ and $R^{n}= [0,1]^{n}$, justifying the 
terminology that $R^{n}\to M$ is a singular {\em cube} in $M$.) 
 
Note that the described ``\pp -formation''  establishes for each $k$  
 a bijection $p_{k}$ between the set $M_{(k)}$ of 
infinitesimal $k$-simplices in $M$, and the set $M_{[k]}$ of 
infinitesimal $k$-dimensional \ppa\ in $M$. So we have the same 
phenomenon as for affine space: a graded set which carries structure 
of (symmetric) simplicial complex, as well as of  cubical complex (with 
reversions and transpositions, and even with subdivisions).

The reader should be warned that the infinitesimal \pg s which arise 
in  this way are more special than the \pg s through which a notion 
of affine connection may be codified, as in \cite{CTC} (1.4); there, 
$x,y,z$ is supposed to satisfy $x\sim y$ and $x\sim z$, but not 
$y\sim z$; without the assumption of $y\sim z$, 
the formation of the \pg\ spanned by $x,y,z$ is an 
added {\em structure}, not canonical.

\section{Higher connections, and their curvature}\label{higher}
It is desirable that  certain mathematical structures  can be encoded 
as {\em morphisms} in a suitable category.
In the context of SDG, a certain general notion of {\em connection} 
can been so 
encoded: recall \cite{CCBI} that if $G$
 is a groupoid with $G_{0}=M$ a manifold, then a {\em connection} in 
$G$ may be construed as a  morphism of reflexive symmetric 
graphs $M_{(1)}\to G$ over $M$. 

We now describe how a certain notion of {\em higher} connection arises in 
a similar way as a  morphism in a suitable category. Note that a 
reflexive symmetric graph is the same as a 1-cubical set with 
reversions. Recall that a cubical groupoid has an underlying cubical 
set with reversions.

\begin{definition}Let $G_{\bullet}$ be a cubical groupoid with object 
set $M$. An {\em $\omega$-connec\-tion} in $G$ is a map $\nabla$ of cubical sets with 
reversion $M_{[\bullet ]}\to G$ over $M$.
\end{definition}
So for each infinitesimal $k$-dimensional \pp\ $P$ in $M$, we have a $k$-cell $\nabla 
_{k}
(P)$ in $G$, in a way which is compatible with face-, degeneracy-, 
and reversion-operators.

\medskip

We have truncated versions, for any finite $k$:

\begin{definition}Let $G$ be a $k$-cubical groupoid with object 
set $M$. A 
{\em $k$-connec\-tion}  is a map $\nabla$ of cubical sets with 
reversion $M_{[[k ]]}\to G$ (preserving $M$) (where $M_{[[k ]]}$ denotes 
the  $k$-truncation of 
$M_{[\bullet ]}$).
\end{definition}
Since $M_{[[k]]}$ is the $k$-truncation of $M_{[[k+1]]}$, we get from 
the adjointness $sk\dashv (-)'$ a bijective correspondence
between the following two kinds of entries

\medskip

\begin{center} \begin{tabular}{ll}
$k$-connections & $\nabla :M_{[[k]]} \to G$\\
\hline
$k+1$-connections &$\widehat{\nabla}  :M_{[[k+1]]}\to G'$.
\end{tabular}
\end{center}

\medskip

\noindent In dimensions $\leq k$, $\nabla$ and $\widehat{\nabla }$ 
agree; in dimension $k+1$, 
$\widehat{\nabla } 
_{k+1}: M_{[k+1]}\to G' _{k+1}$ is what we shall call the {\em formal 
curvature} of $\nabla$. Or, we may call  $\widehat{\nabla }$ itself the 
formal curvature of $\nabla$. With this twist of terminology, the 
formal curvature of a $k$-connection is a $k+1$-connection, and this 
process may then be iterated; this is exploited in the formulation if 
a ``Formal Bianchi Identity'' below, Theorem \ref{bianchi}.

(The view that the curvature of a $k$-connection is a flat 
$k+1$-connection has also been developed by U.\ Schreiber 
\cite{schreiber}, in the context of connections as ``transport along 
paths''.)

\medskip

 Thus, for $k=1$, this means that
for $x\sim y$ in $M$, $\nabla (P(x;y))=\nabla (x,y)$ is an arrow 
$x\to y$ in $G$, and for an infinitesimal \pg\ $P=P(x;y,z)$ with 
fourth vertex $u$, $\widehat{\nabla }(P)$ is the square (shell) in $G$
\begin{equation}\begin{diagram}
z&\rTo ^{\nabla (z,u)}&u\\
\uTo^{\nabla (x,z)}&&\uTo _{\nabla (y,u)}\\
x&\rTo_{\nabla (x,y)}&y 
\end{diagram}\label{aaa}\end{equation}
and this diagram  is only commutative for all such $P$ if $\nabla$ is 
curvature free (flat), in the sense of \cite{CCBI}.
 Thus $\widehat{\nabla} $ in this case encodes the 
information of the curvature of $\nabla$, which is the reason for the 
name ``formal curvature''. In this case ($k=1$), we may get the 
standard \footnote{it is not really ``standard''; but for the 
corresponding simplicial notion of connections, as in \cite{CCBI}, 
 curvature was defined in terms of a similar cyclic composite $R(P):=\nabla (x,y).\nabla 
(y,z)\nabla (z,x)$. A four-fold cyclic composite for defining 
curvature cubically has also been considered by \cite{nish}.} 
(synthetic)
curvature as the function $R$ which assigns to $P$ the arrow $u\to u$ 
given by taking the cyclic composite of the arrows in the displayed 
diagram,
\begin{equation}R(P)=\nabla (u,y).\nabla(y,x).\nabla (x,z).\nabla 
(z,u).\label{cyclic}\end{equation} 
For higher 
$k$, one does not have an analogous construction without a further 
structure on the cubical groupoid $G_{\bullet}$; the structure needed 
is  BSH-connections. Any 1-groupoid carries  canonical such (because 
it is here a void concept); if $G$ 
is an $n$-groupoid carrying  BSH-connections, then $\widehat{G}$ carries
them, too, in a canonical way, see \cite{BH}. 

Recall from \cite{BH} that for a cubical groupoid $G$ with BSH connections, we have 
``folding operations'' ${\bf \phi}:G_{k}\to Cr_{k}(G)$; $Cr _{k}(G)$ 
being, for 
$k\geq 2$, a certain family of groups indexed by $G_{0}$. 
We can use these to reorganize the notion of curvature: 
Let $G$ be a $k$-cubical groupoid with $G_{0}=M$ a manifold, and let 
$\nabla$ be a $k$-connection in it. We  describe the ``real'' 
curvature of it in terms of the formal curvature $\widehat{\nabla}: 
M_{[[k+1]]}\to G'$ previously described, by composing the formal 
curvature with the folding: 

\begin{definition}The {\em } curvature $R_{\nabla}$ of $\nabla$ is the composite map
$$\begin{diagram}[midshaft]M_{[k+1]}&\rTo ^{\widehat{\nabla}_{k+1}}&G' 
_{k+1}&\rTo ^{{\bf \phi}}&Cr_{k+1}(G').
\end{diagram}$$ The $k$-connection $\nabla$ is called {\em flat} if 
its curvature is 0.
\end{definition}
The reader may prefer to further compose with the map $\delta : Cr 
_{k+1}(G')\to Cr_{k}(G')=Cr_{k}(G)$, which is monic, so does not lose 
any information.

\medskip

For $k=1$: consider a 1-connection in an ordinary groupoid $G$. Then 
 $G'$ has a canonical 
BSH-connection, as described
in Section \ref{BSHT}. The corresponding folding operation ${\bf \phi}: 
G'_{2}\to Cr_{2}(G') $ associates to a square in $G$ the cyclic 
composite of its four constituents, as in (\ref{cyclic}). 
($Cr_{2}(G')$ is the $M$-indexed family of the vertex groups of $G$). 
Thus in particular, the connection is flat iff for any infinitesimal \pg , 
the square exhibited in (\ref{aaa}) 
commutes, which is the (cubical rendering of the simplicial) flatness 
notion as described in \cite{CCBI}.

\medskip

Consider again  is a $k$-cubical groupoid $G$ with BSH-connections.
 Then $G'$ is a $k+1$-cubical 
groupoid with BSH-connections, and $G''$ is a k+2-cubical groupoid 
with BSH-connections. 

\medskip

We 
have
\begin{sloppypar}
\begin{thm}[Formal Bianchi Identity] If $\nabla$ is a $k$-connection in 
a $k$-cubical group\-oid $G$ with BSH-connections, 
then the $k+1$-connection $\widehat{\nabla}$  in the 
$k+1$-cubical groupoid $G'$ 
(the formal 
curvature of $\nabla$) is flat.
\label{bianchi}\end{thm}
\end{sloppypar} 
\noindent {\bf Proof.} This is an immediate consequence of the fact that 
$Cr_{k+2}(G'')$ consists of trivial groups $\{ 0 \}$, by Proposition 
\ref{BH5.4}. 

\medskip

\medskip

We would like to describe what this means for $k=1$, entirely in terms 
of $G$, without 
reference to the derived structures $G'$, $G''$, not to speak of the 
folding operations ${\bf \phi}$, which are rather complicated. It is, 
for a 3-dimensional \pp\ $P$ in $M$, entirely a question of the values of 
$\nabla$ on the 12 1-dimensional edges of $P$; they form the 12 edges 
of a cube-shaped diagram ${\bf x}$ in $G$. Now, Proposition \ref{BH5.4} for $n=1$ 
can be read as an identity which holds for all such cube-shaped 
diagrams in all 
groupoids. The Homotopy Addition Lemma (\cite{BH} Lemma 7.1) 
provides simpler expressions for some of the values of the folding 
operations ${\bf \phi}$; for the present case, the value of $\phi$ 
given in loc.cit.\ (Lemma 7.1, case $n=2$): it is the right hand side of 
the 
following 30-letter identity for 12 elements and 
their inverses in a 
groupoid; this identity holds whenever the book-keeping conditions make the 
expression meaningful. We need to give names to the 12 elements 
involved i.e.\ to the edges of the cube-shaped diagram. We number the vertices of 
the cube by 
the integers 0, \ldots ,7; the vertex corresponding to an integer $0, 
\ldots ,7$ is 
then the one whose coordinates in the cube are the digits of the 
integer when written (with three digits) in binary notation, thus 
e.g.\ 3 corresponds to 
the vertex $(0,1,1)$ since 3 in binary notation is 011. Then the 
arrows (edges) are denoted by their domain and codomain; e.g. $37$ 
denotes the arrow from vertex 3 to vertex 7, and $73$ denotes its 
inverse. Note that for instance the symbol $06$ does not correspond 
to an edge of the cube: the twelve edges in question are
01,02,04,13,15,23,26,37,45,46,57,67.
The identity then is  the following; the parentheses can be ignored, but are placed 
there so that one can see each of the six curvatures, corresponding to the 
six faces of the cube. The remaining six arrows perform the needed three 
conjugations which bring those three curvatures 
that live  at other vertices than 7 over 
to vertex 7:

\begin{equation}\begin{split}
\id _{7}=\; &\mbox{ (76 64 45 57) 75 (54 40 01 15) 57 (75 51 13 37) } \\
&\mbox{ \quad \quad 73 (31 10 02 23) 37 
(73 32 26 67) 76 (62 20 04 46) 67}
\end{split}\label{cube1}\end{equation}

It is analogous to the Ph.\ Hall's ``14 letter identity'' used in 
\cite{CCBI} for a simplicial-combinatorial proof of the Bianchi 
identity. A variant of the above identity from \cite{BH}  was found independently bi Nishimura 
\cite{nish}, who also used it to give a proof of the Bianchi identity 
for connections in groupoids.

For  future reference, we rewrite the formula (\ref{cube1}) 
in a (hopefully self-explanatory) notation, with $R$ denoting 
curvature, and exponents denoting conjugation:
\begin{equation}
\id _{7}=R(456).R(041)^{57}.R(153).R(012)^{37}.R(236).R(024)^{67};
\label{cubeR}\end{equation}
here, for instance $(024)$ denotes the square with vertices 0,2,4 and 
6, and 
similarly for the other five $R$-expressions.

\section{Holonomy; Stokes' Theorem}\label{holonomy}Let $M$ be a manifold, and let $G$ be an 
$n$-cubical groupoid with $G_{0}=M$, so that we may talk about 
$n$-connections in $G$. We intend to look for ``integral'' 
 aspects of such connections (i.e.\ {\em global} aspects, rather than infinitesimal).

Recall that we have the cubical set $S_{[\bullet]}(M)$ of singular 
cubes in $M$, and we have its 
$n$-truncation $S_{[[n]]}(M)$.  Its set of 0-cells may be identified with $M$ 
itself. We may therefore consider morphisms of $n$-cubical sets
 $$\Sigma : S_{[[n]]}(M)\to G $$
with $\Sigma _{0} = \id _{M}$.
Recall the notion of subdivisions of 
singular cubes, cf.\ Section \ref{subdivision}. 
\begin{definition}An {\em abstract holonomy} in $G$ is a morphism 
$\Sigma$ of $n$-cubical sets $S_{[[n]]}(M)\to G$ which furthermore has the 
following subdivision property: if an $f\in S_{k}(M)$ subdivides in 
the $i$th direction into 
$f'$ and $f''$, then $$\Sigma (f)= \Sigma (f') +_{i}\Sigma (f'').$$
\end{definition}

 The subdivision property implies that an abstract holonomy $\Sigma$ preserves the 
reversion structure present in both $S_{[[n]]}$ and in $G$.

Recall that the $n$-cubical set $S_{[[n]]}(M)$ also carries the 
structure of transpositions $\tau _{i}$; suppose that the $n$-cubical 
groupoid is equipped with such structure as well. Then we may ask that
an abstract holonomy preserves the transpositions. We then have a 
stronger notion:
\begin{definition}An abstract holonomy $\Sigma$, as in the previous 
definition, is called {\em alternating} if it preserves the 
transposition structure.
\end{definition}
In Section \ref{connections} below, we shall prove that a differential $n$-form 
$\omega $ on 
$M$ gives rise to an $n$-connection (in a certain $n$-cubical 
groupoid), and that integration of $\omega$ over singular $n$-cubes 
in $M$ is an alternating holonomy.

\medskip
The $n$-cubical set $M_{[[n]]}$ of infinitesimal \ppa\ in 
$M$ of dimension $\leq n$ embeds into the $n$-cubical set 
$S_{[[n]]}(M)$, (\ref{in}): $i_{n}: M_{[[n]]}\to S_{[[n]]}(M)$. 
It is clear that the 
truncation functor $tr: Cub _{n+1}\to Cub_{n}$ takes $i_{n+1} $ to 
$i_{n}$. Therefore, the following diagram is commutative, by virtue 
of naturality of the hom-set bijections for the adjointness $tr 
\dashv (-)'$. 
\begin{equation}\begin{diagram}\hom 
(S_{[[n+1]]}(M), G')&\rTo ^{i_{n+1}^{*}}&\hom (M_{[[n+1]]}, G')\\
\uTo ^{\widehat{}}_{\cong}&&\uTo ^{\cong}_{\widehat{}}\\
\hom
(S_{[[n]]}(M), G)&\rTo _{i_{n}^{*}}&\hom (M_{[[n]]}, G),
\end{diagram} \label{pre-stokes}\end{equation} where the vertical 
maps are the hom set bijections $(-)^{\widehat{}}$.

The $\hom$s in this diagram are the hom-set formation for the 
category $Cub_{n+1}$ in the upper line, and $Cub_{n}$ in the lower 
line; However, we    have a similar diagram
if $\hom$ denotes the hom sets for the categories of 
cubical sets with reversion, or cubical sets with reversion and  
transpositions (provided, in the latter case, that $G$ is 
equipped with transpositions). For, the adjointness lifts to these 
more structured categories, and the inclusion $i_{n}: M_{[[n]]}\to 
S_{[[n]]}(M)$  preserves this structure 
(canonically present in the $n$-cubical sets $M_{[[n]]}$ and  
$S_{[[n]]}(M)$);
and similarly for $i_{n+1}$.
\medskip

If  $\nabla : M_{[[n]]}\to G$ is an $n$-connection, one may ask 
whether an abstract holonomy $\Sigma :S_{[[n]]}(M)\to G$ extends $\nabla$, 
in other words, one may ask whether
$$\nabla (P(x_{0}; x_{1}, \ldots ,x_{k}))=\Sigma ([[x_{0}, \ldots 
,x_{k}]])$$
for any infinitesimal $k$-simplex $(x_{0}, \ldots ,x_{k})$ in $M$ 
($k\leq n$). Since an abstract holonomy
by definition has the subdivision property, a necessary condition that $\nabla$ 
extends into an abstract holonomy is that $\nabla$ itself has 
the subdivision property; I don't know whether this in automatic in 
general, but I know cases where it is, see Section \ref{integration} below.

We say that $\nabla$ is {\em integrable} if there is a 
{\em unique} abstract holonomy $\Sigma$ extending $\nabla$; in this case, 
the $\Sigma$ deserves the name {\em the holonomy of $\nabla$} 
and is denoted $\int \nabla$; its 
value on a singular $k$-cube $f$ is denoted by $\int _{f}\nabla$.

If the $n$-cubical groupoid $G$ is furthermore equipped with 
transposition structure, we have a variant of the notion of 
integrability: namely we say that $\nabla$ is {\em a-integrable} if 
there is a unique {\em alternating} holonomy $\Sigma$ extending 
$\nabla$, in which case again this $\Sigma$ is denoted $\int \nabla$. 

\medskip

Examples of $n$-cubical groupoids $G$, in which every connection is 
a-integrable, are given in Section \ref{integration} below; examples of 
1-cubical groupoids (=groupoids), where every connection is 
integrable, were given in \cite{CPCG}.

\medskip
Let $G$ be an $n$-cubical groupoid, and $\Sigma : S_{[[k]]}(M)\to G$ a 
morphism of $n$-cubical sets. 
\begin{lemma} If $\Sigma$ has the 
subdivision property, then so does $\widehat{\Sigma}: S_{[[k+1]]}(M)\to 
G'$.\label{su}\end{lemma}
{\bf Proof.} It suffices to check dimension $n+1$. Let $f\in 
S_{[n+1]}(M)$ $(i,s)$-subdivide into $f'$, $f''$. We must prove
$$\widehat{\Sigma }(f)= 
\widehat{\Sigma}(f')+_{i}\widehat{\Sigma}(f').$$
As elements of $G'_{n+1}$, both sides are determined by their faces, 
so it suffices to see for each $j=1, \ldots ,n+1$ and $\alpha =0,1$ 
that
$$\partial ^{\alpha}_{j}(\widehat{\Sigma}(f))=\partial ^{\alpha}_{j}
\bigl( \widehat{\Sigma}(f')+_{i}\widehat{\Sigma}(f'')\bigr) .$$
There are three cases, $j<i$, $j>i$ and $j=i$. Consider the first 
case. 

We have, by construction of $\widehat{\Sigma}$,
$$\partial ^{\alpha}_{j}(\widehat{\Sigma}(f))=\Sigma \partial 
^{\alpha}_{j}(f).$$
Now by Proposition \ref{subdivisionp}, $\partial 
^{\alpha}_{j}(f)$ $(i-1,s)$-subdivides into $\partial 
^{\alpha}_{j}(f')$ and $\partial 
^{\alpha}_{j}(f'')$, and since $\Sigma$ has the subdivision property, 
we get for the right hand side of the above equation the expression
$$\Sigma (\partial 
^{\alpha}_{j}(f')) +_{i-1} \Sigma (\partial 
^{\alpha}_{j}(f'')) = \partial ^{\alpha}_{j}(\widehat{\Sigma}(f')) 
+_{i-1}\partial ^{\alpha}_{j}(\widehat{\Sigma}(f'')) =\partial 
^{\alpha}_{j}\bigl( (\widehat{\Sigma}(f')+_{i} 
\widehat{\Sigma}(f'')\bigr) ,$$
the last equality sign by one of the rules for cubical groupoids  
(\cite{BH} (1.3.i)).

This proves the desired equation for the case 
$j<i$. The case $j>i$ is similar, except that $i$ does not change to $i-1$ 
during the process. Finally, consider $i=j$. It divides into two 
cases, $\alpha =0$ and $\alpha =1$. For $\alpha =0$, we have
$$\partial ^{0}_{i}\bigl( 
\widehat{\Sigma}(f')+_{i}\widehat{\Sigma}(f'')\bigr) =
\partial ^{0}_{i}\bigl( \widehat{\Sigma}(f')\bigr) = \Sigma 
\bigl( \partial^{0}_{i}(f')\bigr) = \Sigma \bigl( 
\partial^{0}_{i}(f)\bigr) ,$$
using (\ref{subd2}) twice. But this equals $\partial 
^{0}_{i}(\widehat{\Sigma}(f))$, and proves the desired equality.
 The case $\alpha =1$ is similar, using 
(\ref{subd3}). 
\medskip  

 Consider now an $n$-cubical 
groupoid $G$, with $G_{0}=M$ a manifold.
\begin{thm}[Formal Stokes' Theorem] Consider an $n$-connection $\nabla : 
M_{[[n]]}\to G$ 
in $G$; if 
it is integrable,
 then 
so is its formal curvature $\widehat{\nabla}: M_{[[n+1]]}\to G'$; and 
$$\int \widehat{\nabla} =(\int \nabla)^{\widehat{}}.$$
If $G$ is provided with transpositions, the same holds if the word 
``integrable'' is replaced by ``a-integrable''.
\label{Stokes}\end{thm}
{\bf Proof.} Let $\nabla : M_{[[n]]}\to G$ be an integrable 
connection, with holo\-nomy $\int \nabla :S_{[[n]]}(M)\to G$. To say that $\int 
\nabla $ extends $\nabla$ is to say that $i_{n}^{*}(\int \nabla ) 
=\nabla$, and from the commutativity of the diagram 
(\ref{pre-stokes}), it follows that $i_{n+1}^{*}(\int \nabla 
)\widehat{ } = 
\widehat{\nabla}$, which is to say that $(\int \nabla 
)\widehat{ }\; $ extends $\widehat{\nabla}$. Also, since $\int 
\nabla$ has the subdivision property, it follows from Lemma \ref{su} that 
so does $(\int \nabla 
)\widehat{ }\; $. This proves the existence of an abstract holonomy 
extending $\widehat{\nabla}$. The uniqueness is essentially trivial: 
elements in $\hom (S_{[[n+1]]},G')$ ($i=1,2$) 
correspond under the hom-set bijection displayed in 
(\ref{pre-stokes}) to $\Sigma _{i}\in \hom(S_{[[n]]}(M),G)$, i.e.\ 
they are of the form $\widehat{\Sigma}$. Given thus $\widehat{\Sigma}_{i} 
\in  \hom (S_{[[n+1]]},G')$ ($i=1,2$), if they both have the 
subdivision property, then so do $\Sigma _{1}$ and $\Sigma _{2}$; and if 
the $\widehat{\Sigma}_{i}$ both restrict to $\widehat{\nabla}$, then both 
$\Sigma_{i} $s restrict to $\nabla$, and hence $\Sigma _{1}= \Sigma 
_{2}$, by the 
uniqueness assumption (integrability of $\nabla$). 
But then also $\widehat{\Sigma}_{1} =\widehat{\Sigma}_{2}$.

The proof of the assertion about a-integrability is similar, using 
that the diagram (\ref{pre-stokes}) also may be read with $\hom$ 
denoting the category of cubical sets with transpositions and 
reversions.
-- This proves the Theorem.

\medskip

The reason for calling it a ``Stokes' Theorem'' is that, for given 
$f\in S_{n+1}(M)$, $\int _{f} \widehat{\nabla}$ may be read as the 
integral of the coboundary (=curvature) $\widehat{\nabla}$ 
of $\nabla$ over the 
$n+1$-cube $f$, whereas $(\int \nabla )\widehat{} (f)$ may be read as 
the integral of $\nabla$ over the boundary shell of $f$, recalling 
that the values of any $\widehat{h}$ on an $n+1$-cube is obtained in 
terms of the values of $h$ on the boundary shell of the cube. This 
will be elaborated in Section \ref{integration}, 
where we shall consider Stokes' 
Theorem for differential forms.

\section{Cubical-combinatorial differential forms}\label{forms}

We shall give a short exposition of combinatorial differential forms, 
both simplicially and cubically. A fuller account may be found in 
\cite{kompendium}.
 
Let $M$ be a manifold. Recall from Section \ref{infinitesimal} that we have for each $n$ a bijection 
$p_{n}$ from the set $M_{(n)}$ of infinitesimal $n$-simplices in $M$
 to the set $M_{[n]}$ of infinitesimal $n$-dimensional \ppa\ in $M$. 
Thus there is also a bijective correspondence $p_{n}^{*}$ between 
functions $M_{[n]}\to R$ and functions $M_{(n)}\to R$. One of the 
 equivalent ways to say  that  such a 
function $M_{(n)}\to R$ is a {\em combinatorial differential $n$-form} 
(or simplicial-combinatorial differential $n$-form on $M$) is that it 
takes value 0 on any $(x_{0}, \ldots ,x_{n})$ where $x_{0}=x_{i}$ for 
some $i>0$. Simplices of this form are precisely those that under 
$p_{n}$ correspond to values of the cubical degeneracy operators 
(such we call {\em degenerate} cubes). 
Therefore we put

\begin{definition}A {\em cubical-combinatorial differential $n$-form} 
(briefly, a cubical $n$-form) on $M$ is a function $M_{[n]}\to R$ 
which take degenerate $n$-cubes to 0.
\end{definition} 
(When we in the following say ``$n$-form'' without further 
decoration, we mean 
``cubical-combinatorial differential $n$-form''.)

 So the bijection $p_{n}^{*}$ 
restricts to a bijection between simplicial-combi\-na\-torial $n$-forms 
and cubical-combinatorial $n$-forms.

Now both $M_{(\bullet )}$ and $M_{[\bullet ]}$ are {\em symmetric} 
simplicial (resp.\ cubical) sets: the 
symmetric group ${\mathfrak S}_{n+1}$ acts on 
$M_{(n)}$, and ${\mathfrak S}_{n}$ acts on $M_{[n]}$. Also 
${\mathfrak S}_{n}$ acts on $M_{(n)}$, as those permutations in 
${\mathfrak S}_{n+1}$ which 
keep the first vertex fixed.
Simplicial-combinatorial $n$-forms are known to be {\em alternating}, 
in the sense that $\omega (\tau (x)) = \sign (\tau)\cdot \omega (x)$ 
for any $\tau \in {\mathfrak S}_{n+1}$, $x\in M_{(n)}$.
  
As we observed in Section \ref{complex} for the case of affine 
simplices and \ppa , $p_{n}$ preserves the action of ${\mathfrak S}_{n}$; this 
also holds for infinitesimal simplices/ \ppa ; so it immediately follows that 
cubical-combina\-to\-rial $n$-forms are likewise ${\mathfrak S}_{n}$-alternating.

These assertions, however, do not exhaust the symmetry properties of 
neither the simplicial $n$-forms, nor the cubical $n$-forms; 
simplicial $n$-forms are not only alternating w.r.to the action of 
${\mathfrak S}_{n}$, but with respect to the action of ${\mathfrak 
S}_{n+1}$ on $M_{(n)}$; and cubical $n$-forms are alternating w.r.to 
the ``reversion'' structure which the set $M_{[n]}$ carries. For this 
latter assertion, let us be more explicit. By a coordinate argument, 
one may prove

\begin{prop}Let $\omega$ be a cubical-combinatorial $k$-form on $M$. 
Then for any infinitesimal \pp\ $P$, $ \omega (-_{i}P)=-\omega (P)$,
where $-_{i}P$ is the reversion of $P$ in the $i$th direction $(i=1, 
\ldots ,k)$.
\label{rever}\end{prop}

We refer to \cite{kompendium} for a proof. It is essentially the same 
as the proof that gives the ``extra'' symmetry property (${\mathfrak 
S}_{n+1}$-alternating property) for simplicial $n$-forms.

\medskip

On $R^{n}$, there is a canonical $n$-form, the ``volume form'' 
$\Vol$. For an infinitesimal \pp\ $P=P(x_{0};x_{1},\ldots ,x_{n})$, 
it is given by the formula
$$\Vol (P):=\det (x_{1}-x_{0}, \ldots ,x_{n}-x_{0}),$$
the determinant of the $n\times n$ matrix whose columns are 
$x_{i}-x_{0}$. (The form $\Vol$ may also be defined as $dt_{1}\wedge \ldots \wedge 
dt_{n}$, with the wedge product as described in \cite{kompendium}.)

We shall  need the following result, see \cite{kompendium}:

 \begin{prop}For every cubical-combinatorial 
$n$-form $\theta$ on $R^{n}$, there exists a 
unique function $\widehat{\theta }:R^{n} \to R$ such that for any infinitesimal \pp\ 
$P=P(x_{0};x_{1}, \ldots ,x_{n})$
$$\theta (P) = \widehat{\theta }(x_{0})\cdot \det (x_{1}-x_{0}, \ldots 
,x_{n}-x_{0}),$$
i.e.\ such that $\theta = \widehat{\theta}\cdot \Vol $.\label{theta}\end{prop}

Since arbitrary functions $f:N\to M$ between manifolds take 
infinitesimal \ppa\ to infinitesimal \ppa\, we have, just as for 
simplicial-combinatorial forms, that the set of cubical $n$-forms on 
manifolds depends contravariantly on the manifold; thus, if $\theta$ 
is a cubical $n$-form on $M$, we get a cubical $n$-form $f^{*}(\theta 
)$ on $N$, with
$$f^{*}(\theta )\bigl (P(y_{0};y_{1},\ldots ,y_{n})\bigr ):= \theta 
\bigl( P(f(y_{0});f(y_{1}), \ldots ,f(y_{n})\bigr) ,$$
where the $y_{i}$s are mutual neighbours in $N$.  We are going to use this for
the case where $f$ is 
a map $R^{n}\to M$ given by some infinitesimal $n$-\pp\ in $M$;
namely, we have the following Proposition from \cite{kompendium}. It is proved by a 
Taylor expansion argument, together with the product rule for 
determinants:

\begin{prop}Let $\omega$  be a cubical $n$-form on a manifold $M$, 
and let $P=P(x_{0};x_{1}, \ldots ,x_{n})$ be an infinitesimal 
$n$-dimensional \pp\ in $M$ . Then we have the following equality of
 $n$-forms on 
$R^{n}$: 
$$[[x_{0}, \ldots ,x_{n}]]^{*}(\omega ) = \omega (P) \cdot \Vol .$$
\label{subst}\end{prop}
Note that the right hand side is a ``constant'' $n$-form on $R^{n}$ 
in the sense that the function $\widehat{\theta}$ corresponding, by 
Proposition \ref{theta}, to it, 
 is a constant function (constant with value $\omega (P)$). 

\medskip

In Section \ref{integration}, we shall consider the question of integrating 
$n$-forms  along singular $n$-cubes.

\medskip

For a simplicial-combinatorial $n$-form $\omega$, one constructs a 
simplicial- combinatorial $n+1$-form $d\omega$, by the usual formula 
for coboundary of simplicial cochains; see \cite{CCBI} \S 4. For 
cubical-combinatorial $n$-forms, one also has a coboundary operator, 
by the ususal formula for coboundary of cubical cochains (as in 
\cite{HW} \S 8.3, say). These two 
coboundaries match modulo a factor $n+1$: denoting the coboundaries 
respecively by $d_{s}$ and $d_{c}$ (for ``simplicial'' and 
``cubical'', respectively), the formula is
$$d_{s}(\omega )= \frac{1}{n+1} d_{c}(\omega ),$$
(where we omit the bijections $p_{n}^{*}$ and $p_{n+1}^{*}$ from 
notation). The cubical formula can be seen, when working in a 
coordinatized situation, to give the standard formula for exterior 
derivative of ``classical'' differential forms (via a well known 
correspondence between combinatorial and classical differential 
forms, see \cite{kompendium}). 

We  have have 
$d(f^{*}(\theta )) = f^{*}(d\theta )$.

\section{Connections as generalized forms}\label{connections}
We shall see how the notion of ``connection in an $n$-cubical groupoid'' 
contains as a special case the notion of differential $n$-form. This 
comes about by considering certain ``constant'' $n$-cubical groupoids.
(For $n=1$, this was described in \cite{CPCG}; here the ``value 
group'' $A$ 
need not be assumed commutative.)

Let $M$ be any set (manifold), and $(A,+)$ any abelian group (ultimately, we shall 
be interested in the case where $A=R$, the number line). For $n\geq 1$, we 
may consider the $n$-cubical groupoid $M_{n}(A)$ whose $k$-cells 
for $k<n$ form the set $M^{2^{k}}$; and the set of $n$-cells is taken 
to be $M^{2^{n}}\times A$. (It is the cartesian product in the 
category of $n$-cubical groupoids of the codiscrete $n$-cubical 
groupoid on $M$, on the one side, and, on the other side, 
the $n$-cubical groupoid with $A$ as set of 
$n$-cells,  a one-point set as set of $k$-cells, for $k<n$, and 
where all 
$n$ compositions $+_{i}$ of $n$-cells are taken to be $+$.

The $n$-cubical groupoid $M_{n}(A)$ has BSH connections (using $0\in 
A$).

The 
category of $n$-cubical groupoids with BSH-connections is, by 
\cite{BH}, equivalent 
to the category of crossed $n$-comlexes over groupoids, as described in 
loc.cit.; such a thing consists in an ``ordinary'' groupoid $G\tto 
M$, and for each $n\geq 2$, an $M$-indexed family of ($n$-truncated) chain complexes
$C_{n}\stackrel{\delta}{\to} C_{n-1}\to \ldots C_{3}\stackrel{\delta}{\to} C_{2}$ 
(where $C_{2}$ is the 
family of vertex groups of $G\tto M$; $C_{2}$ need not 
 be abelian; the groupoid $G\tto M$ acts on these chain 
complexes in a compatible way. We shall not need a more complete 
description, but just describe the crossed $n$-complex (respectively 
crossed $n+1$-complex) corresponding to $M_{n}(A)$ (respectively 
to $(M_{n}(A))'$). In both cases, the groupoid is the codiscrete 
$M\times M \tto M$, and it acts trivially. For each $m\in M$, the 
corresponding chain complexes are, respectively, $A$ concentrated in dimension $n$ 
(and $0$s elsewhere); and $A$ in dimension $n+1$ and also in dimension 
$n$, with the identity map as $\delta$, and $0$s elsewhere.

For any $n$-cubical groupoid $G$, let us denote the corresponding 
crossed complex by $Cr (G)$. It agrees with $G$ in dimension 0 and 
1. Its $k$-dimensional part $Cr _{k}(G)$ ($k\geq 2$) is a subset of $G_{k}$, 
consisting of those $k$-cells all of whose faces except $\partial 
^{0}_{1}$ are ``totally degenerate'', i.e.\ comes about from a 0-cell 
by applying $k$ degeneracies. Then (the restriction of) $\partial 
^{0}_{1}$ will serve as $\delta : Cr _{k}(G)\to Cr_{k-1}(G)$. 
According to \cite{BH}, there is 
for each $k$ a retraction ${\bf \phi} : G_{k}\to Cr _{k}(G)$; it is a 
somewhat complicated ``folding operation''. However, for
$M_{n}(A)$ and $(M_{n}(A))'$, the folding operations
 are easy to describe. In dimensions $k<n$, it is the map 
$l:M^{2^{k}}\to M$ which pick out the last vertex of a $2^{k}$-tuple.
In dimension $n$, it is the map $l\times \id _{A}:M^{2^{n}}\times 
A\to M\times A$ (these descriptions apply both to $M_{n}(A)$ and 
to $(M_{n}(A))')$, which agree in dimensions $\leq n$). Finally, 
for $(M_{n}(A))'$, we need to describe $\phi$ in dimension $n+1$. 
It follows from the ``Homotopy Addition Lemma'', \cite{BH} Lemma 7.1, 
that $\phi$ here may be described as follows. An element in
$(M_{n}(A))'_{n+1}$ consists of a $2^{n+1}$-tuple ${\bf x}$ of elements 
from $M$, together with a $2(n+1)$ tuple of elements from $A$, one 
for each face operator $\partial ^{\alpha}_{i}$ ($\alpha = 0,1$, 
$i=1, \ldots ,n+1$). Denoting the element $\in A$ corresponding to
$\partial ^{\alpha}_{i}$ by 
$a^{\alpha}_{i}\in A$, $\phi$ gives the value 
\begin{equation}\bigl( l(x), \sum _{i=1}^{n+1} (-1)^{i}\{ 
a^{1}_{i}-a^{0}_{i}\} \bigr) .\label{hal}\end{equation}
Note that in the quoted Lemma, it is $\delta \circ \phi $ rather than 
$\phi$ which is 
described, but in our case $\delta$ is an identity map. Also note 
that all the groupoid actions  appearing in the formula, 
can be ignored; the action is trivial. Finally note that the $l(x)\in M$ 
components in the formulae for ${\bf \phi}$ contain no information, all information 
resides in the component $\in A$. 

\medskip
We now consider the case $A=R$, for concreteness; 
for the case $n=1$, non-commutative Lie 
groups may be used for $A$, see \cite{CPCG}). For $n=2$, see the 
Remark at the end of the present Section.
 Then we have
\begin{prop}\label{seven}There is a bijective correspondence between
connections $\nabla$ in the $n$-cubical groupoid $M_{n}(R)$, and differential 
($R$-valued) $n$-forms $\omega$ on $M$. To $\nabla$, the corresponding 
$\omega$ is given as the composite
$$\begin{diagram}[midshaft]M_{[n]}&\rTo ^{\nabla _{n}}&M^{2^{n}}\times R&\rTo 
^{\proj}&R.\end{diagram}$$
Also,  connections in $M_{n}(A)$ are automatically alternating. 
\end{prop}
{\bf Proof.} To see that the exhibited map is a 
(cubical-combinatorial) $n$-form, it suffices to see that it takes 
value 0 on degenerate $n$-cells; but $\nabla$ takes degenerate 
$n$-cells into identity $n$-arrows in the groupoid, and identity 
arrows have 
their $A$-component equal 0. -- On the other hand, given an $n$-form 
$\omega$, we define $\nabla :M_{[[n]]}\to M_{n}(R)$ as follows: in 
dimensions $k<n$, it is just the inclusion map $i: M_{[k]}\to M^{2^{k}}$ 
associating to an infinitesimal $k$-dimensional \pp\ its $2^{k}$-tuple of vertices; 
and in dimension $n$, it associates to an infinitesimal 
$n$-dimensional \pp\ $P$ the $n$-cell given by $(i(P)), \omega (P))\in 
M^{2^{n}}\times R$. This clearly describes a morphism of cubical sets 
with reversions 
It is clear that the correspondence is bijective. Since differential 
$n$-forms $\omega $ 
are known to be automatically alternating, the last assertion of the 
Proposition follows.

\medskip

With $\nabla$ and $\omega$ as in this Proposition, we shall next 
relate the formal curvature $\widehat{\nabla}$ of 
$\nabla$ with the exterior derivative of the  $\omega$.

\begin{prop}\label{eight}With $\nabla$ and $\omega$ as in the previous 
Proposition,  the composite
$$\begin{diagram}[midshaft]M_{[n+1]}&\rTo ^{\widehat{\nabla}_{n+1}}&(M_{n}(R))'_{n+1}&\rTo 
^{{\bf \phi}}&M\times R &\rTo ^{\proj}&R\end{diagram}$$
equals $d\omega $.
\end{prop}
{\bf Proof.} Recall that the cubical-combinatorial $n+1$-form 
$d\omega$ is given by the standard cubical coboundary of $\omega$, 
thus for $P$ an infinitesimal $n+1$-dimensional \pp\,
$$d\omega (P)= \sum _{i=1}^{n+1} (-1)^{i}\{ \omega (\partial 
^{1}_{i}(P)) - \omega (\partial ^{0}_{i}(P))\} .$$
Now the result follows by comparison with (\ref{hal}). 

\medskip

\noindent {\bf Remark.} We finish by a remark on 2-forms with 
non-commutative values. Recall that a crossed complex need not be commutative in dimension 2. 
By the equivalence between cubical groupoids with BSH connections and 
crossed complexes, we have the possibility of deriving a notion of 
2-form with non-commutative values, and its coboundary, by passing 
via the curvature of the corresponding 2-connection. We shall just 
sketch this. (The reader may want to specialize to the case where the 
crossed complex is constant, as we did when we viewed  
$R$-valued $n$-forms as $n$-connections.)

So consider a crossed module (= 2-crossed complex) corresponding to a 
2-cubical groupoid $G$ with BSH-connections. The 0-dimensional part 
is supposed to be a manifold $M$. 
To fix notation, we exhibit this crossed module $Cr(G)$:
$$\begin{diagram}C_{2}&\rTo ^{\delta }&G_{1}&\pile{\rTo \\ \rTo }& M.
\end{diagram}$$
A 2-connection $\nabla $ in $G$ 
gives, via the folding operation $\phi$, rise to the following data with values in the crossed module 
$Cr(G)$:

1) A 1-connection $\nabla _{1}$ in the 1-groupoid $G_{1}\tto M$;

2) A function $\omega$ on $M_{[2]}$ with values in the group bundle $C_{2}$ 
over $M$ 

satisfying, for all infinitesimal \pg s $P(x_{0};x_{1},x_{2})$

$$\delta (\omega (P(x_{0};x_{1},x_{2})) = R(P(x_{0};x_{1},x_{2})).$$
Here, $\omega$ is a group-bundle-valued 2-form, in the sense that 
$\omega (P(x_{0};x_{1},x_{2}))\in C_{2}(x_{3})$, (= the group 
corresponding to the fourth vertex $x_{3}\in M$ of the \pg\ 
$P(x_{0};x_{1},x_{2})$ \footnote{conventions  are a little clumsy 
here -- we 
might have chosen that $\omega (P(x_{0};x_{1},x_{2})\in C_{2}(x_{0})
$, say, but we want to stick to the conventions involved in the 
construction of the folding operations from \cite{BH}.}), and so that 
the value is $0$ if the \pg\ is degenerate (we are using 
additive notation in the groups that make up $C_{2}$). Also, $R$ 
denotes the (real) curvature of the connection $\nabla _{1}$, arising 
from the formal curvature $\widehat{\nabla _{1}}$ by the folding.

Now the formal curvature $\widehat{\nabla}$ of the 2-connection is a 
3-connection in $G'$, which is a 3-cubical groupoid with 
BSH-connections. To this 3-cubical groupoid corresponds the 3-crossed 
complex
$$\begin{diagram}\mbox{Ker}(\delta )&\rInto & C_{2}&\rTo ^{\delta}&G&
\pile{\rTo \\ \rTo }& M.
\end{diagram}$$
(Here, $\mbox{Ker}(\delta )$ is contained in the center of $C_{2}$, 
and it is in 
particular commutative.) So the formal curvature $\widehat{\nabla}$ 
gives by folding rise to a map from $M_{[[3]]}$ to this crossed 
complex, and the only added information in this is the 3-dimensional 
part $M_{[3]}\to \mbox{Ker}(\delta )$; this map, we consider as the 
coboundary of the pair $\omega ,\nabla _{1}$. 
 The explicit formula for this coboundary can be read out of the 
formula (\ref{cubeR}); on a $P\in M_{[3]}$, say 
$P=P(x_{0};x_{1},x_{2},x_{4})$ with vertices $x_{0}, \ldots ,x_{7}$,
$R(456)$ denotes $R(P(x_{4};x_{5},x_{6}))$ etc., and the exponent 
formation, say the exponent $57$ in the second term, denotes the 
action of $\nabla _{1}(x_{5},x_{7}) \in G$ on $C_{2}$.

\section{Integration of differential forms}\label{integration}
We present a synthetic theory of integrals which avoids the use of Riemann sums or 
other approximation techniques and entirely depends on the Fundamental 
Theorem of Calculus. It also, for simplicity, does not depend on any preorder $\leq$ 
on the number line $R$. It depends on one integration axiom, namely: 
{\em to any $f:R\to R$, there exists an $F:R\to R$ with $F' =f$; and such 
$F$ is unique up to a constant}. (The function $F$ is called a {\em 
primitive} of $f$.) With this axiom, one defines one-dimensional 
integrals: $\int _{a }^{b} f := F(b)-F(a)$, for any $a,b$ in $R$. We 
can then also define {\em iterated} integrals: for any $f:R^{2}\to R$, one 
defines $\int _{a_{1}}^{b_{1}} \int _{a_{2}} ^{b_{2}} f$, for any 
$a_{1}, b_{1}$ and $a_{2}, b_{2}$ in $R$, and similarly for $n$-fold 
iterated integrals.

For simplicity of exposition and geometry, 
we shall here restrict ourselves to to the case $n=2$.

Let now $M$ is a manifold and $\omega$ a (combinatorial-cubical) 
2-form on $M$. Then we can define a function 
$$\int \omega : S_{[2]}(M) \to R$$ as follows. Let $f: R^{2}\to M$ be 
a singular square (= singular 2-cube). By Proposition 
\ref{theta}, the 2-form $f^{*}(\omega )$ on $R^{2}$ is of the form
$$f^{*}(\omega ) = \widehat{\theta} \cdot \Vol$$ for a unique 
function $\widehat {\theta}:R^{2}\to R$, and we define
$$\int _{f}\omega := \int _{0}^{1} \int _{0}^{1}\widehat{\theta}.$$
(We also write the iterated integral in this formula as $\iint 
_{I\times I}\widehat{\theta}$.)
Let us temporarily denote the  function $\int \omega : S_{[2]}(M) \to R$
thus defined by the symbol $\Omega$. Then $\Omega$ has the properties (i) and (ii) 
below, which we 
shall use as definition of the notion of {\em abstract surface 
integral}\footnote{The condition (i) is essentially the same as the 
defining property of an ``observable'', as considered in \cite{MR}, 
Definition 16.}:

\medskip

(i) $\Omega (f) = \Omega (f')+ \Omega (f'')$ whenever $f$ subdivides 
into $f'$ and $f''$;

(ii) $\Omega$ is alternating in the sense that $\Omega (f\circ \tau ) 
= - \Omega (f)$

\medskip

\noindent where $\tau :R^{2}\to R^{2}$ interchanges the two coordinates.
All this is as one would expect; it is standard multivariable 
calculus. The following, however, has no place in the standard 
treatment.
\begin{prop}Let $\omega$ be a  cubical-combinatorial $2$-form on a 
manifold $M$, and let $P=P(x;y,z)$ be an infinitesimal \pg . Then
\begin{equation}\int _{[[x,y,z]]}\omega = \omega 
(P).\label{inf-beh2}\end{equation}
\label{IB}\end{prop}
{\bf Proof.}  This follows immediately from Proposition \ref{subst}, since 
$\iint_{I\times I} 
c$ for $c$ a constant equals $c$.

\medskip

As a Corollary, we have that (cubical-combinatorial) 2-forms $\omega$ have the 
subdivision property: if an infinitesimal \pg\ $P$ subdivides into 
$P'$ and $P''$, then $\omega (P)= \omega (P')+\omega (P'')$. For, 
$\int \omega$ has this property, for all \pg s, and for infinitesimal 
\pg s, $\int \omega$ agrees with $\omega$ itself, by (\ref{inf-beh2}).

\medskip

For an abstract surface integral $\Omega$ on $M$, we may ask whether it 
extends a given 2-form $\omega$ on $M$, in the sense that
$$\Omega ([[x,y,z]])=\omega (P(x;y,z))$$
for any infinitesimal \pg\ $P(x;y,z)$. The ``concrete'' surface 
integral $\int \omega$ does extend $\omega$ in this sense, by  
Proposition \ref{IB}. We have more completely:

\begin{thm}For every 2-form $\omega$, $\int \omega$ is the only 
abstract surface integral extending $\omega$. Thus, there is a 
bijective correspondence between abstract surface integrals, and 
2-forms.
\label{ccc}\end{thm} 
{\bf Proof.}  By 
subtraction, this uniqueness assertion is equivalent to the following:
\begin{lemma}If an abstract surface integral $\Phi$ on $M$ has the 
property that it vanishes on all infinitesimal $\pg s$, then it 
vanishes.
\label{van11}\end{lemma}
(Recall that an infinitesimal \pg\ $P(x_{0};x_{1},x_{2})$ may be 
identified with a certain map $[[x_{0},x_{1},x_{2}]]:R^{2}\to M$.)

\medskip

\noindent {\bf Proof.} We first note that if $f:N\to M$ is any map 
between manifolds, an abstract surface integral $\Phi$ on $M$ gives 
rise in an evident way to a an abstract surface integral $f^{*}(\Phi 
)$ on $N$, $f^{*}(\Phi )(g) := \Phi (f\circ g)$ for any $g:R^{2} \to 
N$. Since $f$ takes infinitesimal \pg s to infinitesimal \pg s, we 
have by the assumption on $\Phi$ that $f^{*}(\Phi )$ vanishes on 
infinitesimal \pg s in $N$.

Now let $f:R^{2}\to M$ be an arbitrary element of $S_{2}(M)$; we have 
to see that $\Phi (f)=0$. Since $\Phi (f)= f^{*}(\Phi )(\id )$, and 
$\id$ (the identity map on $R^{2}$) is a \pg\ in the affine space 
$R^{2}$ (namely $P(0;e_{1},e_{2})$, it 
suffices to see that $f^{*}(\Phi )$ vanishes on all \pg s in $R^{2}$.
In other words, the Lemma follows by taking $\Psi = f^{*}(\Phi )$ in 
the following

\begin{lemma}If an abstract surface integral $\Psi$ on $R^{2}$ 
vanishes on all infinitesimal \pg s, it vanishes on all \pg s.
\label{van22}\end{lemma} 
{\bf Proof.} This will follow if we can prove each of the following 
three assertions:

\medskip

\noindent 1) If $\Psi$ vanishes on all infinitesimal \pg 
s, it vanishes on all \pg s with infinitesimal sides.

\noindent (This means \pg s $P(x;y,z)$ with $x\sim y$ and $x\sim z$, but not 
necessarily with $x\sim z$.)

\medskip

\noindent 2) If $\Psi$ vanishes on all \pg s with infinitesimal sides, it 
vanishes on all \pg s which have (at least) {\em one} side infinitesimal.

\medskip

\noindent 3) If $\Psi$ vanishes on all \pg s which have one side infinitesimal, 
it vanishes on all \pg s.

\medskip

The proof of 1) is a piece of infinitesimal algebra from \cite{SDG}, (and does not 
depend on the subdivision property): for fixed $x$, consider the 
function $R^{2}\times R^{2}\to R$ given by $(u,v)\mapsto \Psi 
(P(x;x+u,x+v))$. It is alternating, because of the alternating 
property of $\Psi$, and it vanishes if one of its arguments is 0.
 Its restriction to $D(2)\times D(2)$ extends therefore to a 
bilinear alternating map $g: R^{2}\times R^{2}\to R$. By assumption, 
$g$ vanishes on $\widetilde{D}(2,2)\subseteq D(2)\times D(2)$. But 
bilinear alternating maps $R^{2}\times R^{2}\to R$ are determined by 
their restriction to $\widetilde{D}(2,2)$, see \cite{SDG} I. 16.

The proofs of 2) and 3) are quite similar to each other.  
Let $f(s):=\Psi (P(x;y(s),z)$. 
Now assume that $\Psi$ 
satisfies the assumption of 2), and consider a \pg\  
$P(x;y,z)$ with the side $(x,z)$ infinitesimal. We use notation 
$y(s)$, $u(s)$ as in the statement of (\ref{subdiv}). Then, for $d\in D$,
\begin{equation}f(s+d) = \Psi (P(x;y(s+d),z)) = \Psi (P(x;y(s),z))+\Psi 
(P(y(s);y(s+d), u(s))),\label{van2}\end{equation}
by the subdivision property of $\Psi$ . The last term here vanishes 
because the \pg\  $P(y(s);y(s+d), u(s))$ has infinitesimal sides 
(note $z-x=u(s)-y(s)$, so $u(s)\sim y(s)$). We conclude $f(s+d)=f(s)$ for all $d\in D$, 
so $f'(s)=0$, for all $s$. Since also $f(0)= 
\Psi (P(x;x,z))=0$ ($P(x;x,z)$ being infinitesimal), we conclude that 
$f$ is constant 0. Hence $f(1)=0$, and 
this assertion is equivalent to $\Psi (P(x;y,z))=0$

If we instead had considered a \pg\ with its {\em first} side 
$(x,y)$ infinitesimal, we reduce to the case treated just by using 
the alternating property $\Psi (P(x;y,z))=-\Psi (P(x;z,y))$.

Finally, assume that $\Psi$ satisfies the assumption in 3). We 
consider an arbitrary \pg\ $P(x;y,z)$, and define $y(s)$, $u(s)$ as 
above, and again define $f(s):=\Psi (P(x;y(s),z)$. We then again have 
the equation (\ref{van2}), and now the last term vanishes by 
assumption, because 
$P(y(s);y(s+d), u(s))$ has one side $(y(s),y(s+d))$ infinitesimal.
So again we conclude $f'(s)=0$. Also $f(0)=0$, since $f(0)$ is $\Psi$ 
applied to a \pg\ with one of its sides 0 (hence in particular, with 
an infinitesimal side). So $f$ is constant 0, hence $f(1)=0$, and 
this assertion is equivalent to $\Psi (P(x;y,z))=0$. This proves 3), 
and by combining 1), 2) and 
3), we get the Lemma.

\medskip

\begin{prop}Let $\Phi$ be an abstract surface integral on a manifold 
$M$. Then there exists a unique 2-form $\omega$ on $M$ such that $\Phi = 
\int \omega$.
\end{prop}
{\bf Proof.} Uniqueness is clear; for an infinitesimal \pg\ 
$P=P(x;y,z)$, we must necessarily put $\omega (P):= \Phi ([[x,y,z]])$.
To see that $\int \omega$ and $\Phi$ agree, it suffices by 
Lemma \ref{van11} to see that they agree on infinitesimal \pg s 
$[[x,y,z]]$, but  $\omega$ was defined by that.

\medskip

Consider the $n$-cubical groupoid $M_{n}(R)$.  It has transpositions.
We have in Section 
\ref{connections} seen that there is a bijective correspondence between alternating connections in 
it, and $n$-forms on $M$ (and every connection is automatically 
alternating). The same recipe which gave this 
correspondence also gives a correspondence between abstract 
alternating holonomies in it, and abstract $n$-dimensional surface 
integrals. 
\begin{prop}The n-cubical groupoid $M_{n}(R)$ admits a-integration, 
i.e.\ every n-connection in it (automatically alternating) extends 
uniquely to an alternating holonomy $S_{[n]}(M)\to M_{n}(R)$. Also, 
if the $n$-form $\omega$ corresponds to the $n$-connection $\nabla$,
\begin{equation}\int \omega =\proj \circ \int \nabla .
\label{int}\end{equation}
\label{admits}\end{prop}
{\bf Proof.} In view of the bijective correspondences  which we, for 
this $n$-cubical groupoid, have between 
(alternating) $n$-connections and $n$-forms, and between alternating 
holonomies and abstract $n$-dimensional surface integrals, this follows,
for $n=2$,  immediately 
from Theorem \ref{ccc}. (The proof in other dimensions is essentially 
the same.)

\medskip

Let us finally analyze how a (restricted) version of Stokes' theorem for 
1-forms follows from the results and constructions presented. We 
consider a 1-form $\omega$ on $M$, and the corresponding connection 
$\nabla$ in $M_{1}(R)$. We need to be precise about the crossed 
complex  $Cr (M_{1}(R)')$ correspondingto the 2-cubical groupoid 
$M_{1}(R)'$; in 
dimensions 1 and 2, it is  just $M\times R$, and the $\delta$ is 
the identity map. We can safely ignore the $M$-component, so we shall 
not consider the folding operations $\phi$ themselves, but their composite with 
their projection $\overline{\phi}:M\times R \to R$.  In dimension 1, 
$\overline{\phi}$ itself is just a projection, (recall that 
$M_{1}(R)$ in dimension 1 is $M^{2}\times R$). In dimension 2, it 
follows from the Homotopy Addition Lemma \cite{BH} Lemma 7.1 (case 
$n=1$) that $\overline {\phi}$ to a shell (square of arrows) 
associates a suitable alternating four-fold sum of the numbers associated to the 
four arrows of the square. We need one further property of 
$\overline{\phi}$, namely that it takes composites $+_{1}$ and 
$+_{2}$ in $M_{1}(R)'$ to sum formation $+$ in $R$; this follows from 
\cite{BH}, (4.9)(i). 

So consider a 1-form $\omega$ on $M$, and let $\nabla$ be the 
connection in $M_{1}(R)$ corresponding to it by Proposition 
\ref{seven}. Then by Proposition \ref{eight}, $d\omega 
=\overline{\phi} \circ \widehat{\nabla}$. Since $\int 
\widehat{\nabla}$ extends $\widehat{\nabla}$, 
$\overline{\phi}\circ \int 
\widehat{\nabla}$ extends $\overline{\phi} \circ \widehat{\nabla}$.
Since $\int 
\widehat{\nabla}$ has the subdivision property, it follows from the 
quoted equation (4.9)(i) that also $\overline{\phi}\circ \int 
\widehat{\nabla}$ has the subdivision property. So it has the properties 
which characterize $\int d\omega$, so that we have the first equality 
sign in
$$\int d\omega = \overline{\phi}\circ \int \widehat{\nabla} = 
\overline{\phi} \circ \widehat{\int \nabla},$$
the second equality by the ``Formal Stokes' Theorem \ref{Stokes}. 
Let $f\in S_{2}(M)$ be a singular 2-cube; by the quoted instance of 
the Homotopy Addition Lemma, applied to the shell $\widehat{\int 
\nabla}(f)$,
$\delta (\overline{\phi}(\widehat{\int \nabla}(f))$ equals $\sum \pm 
\overline{\phi}(\int _{\partial^{\alpha}_{i}(f)}\nabla )$, (standard fourfold 
sum) which is 
$\sum \pm \int _{\partial^{\alpha}_{i}(f)}\omega = \omega (\partial 
f)$.

\end{document}